\documentclass{article}
\usepackage[utf8]{inputenc}
\usepackage[margin=.9in]{geometry}
\usepackage{amsfonts,amsmath}
\usepackage{cite}
\def\eq#1{\begin{equation}#1\end{equation}}
\newcommand{\R}{{\rm I\!R}}
\def\rep#1{(\ref{#1})}
\def\scr#1{{\cal #1}}
\newtheorem{theorem}{Theorem}
\newtheorem{lemma}{Lemma}
\newtheorem{proposition}{Proposition}

\newtheorem{assumption}{Assumption}
\newtheorem{remark}{Remark}
\def\qed{ \rule{.08in}{.08in}}

\newcommand{\bbb}{\mathbb}

\newcommand{\1}{\mathbf{1}}
\newcommand{\0}{\mathbf{0}}

\newcommand{\D}{\mathcal{D}}

\newtheorem{definition}{Definition}

\title{On the Analysis of a Continuous-Time Bi-Virus Model
}
\author{Ji Liu, Philip E. Par\'{e}, 
Angelia Nedi\'{c}, Choon Yik Tang, 
Carolyn L. Beck, and Tamer Ba\c sar\thanks{
Ji Liu, Philip E. Par\'{e}, 
Angelia Nedi\'{c}, 
Carolyn L. Beck, and Tamer Ba\c sar are with the Coordinated Science Laboratory at the University of Illinois at Urbana-Champaign.  Choon Yik Tang is with the ECE Department at the University of Oklahoma.  This material is based on research partially sponsored by the National Science Foundation grants ECCS 15-09302, CCF 11-11342, and DMS 13-12907 and ONR Basic Research grant Navy N00014-12-1-0998.  All material in this paper represents the position of the authors and not necessarily that of NSF.
A version of this paper was submitted to the 2016 IEEE Conference on Decision and Control.
}}
\date{}

\begin{document}

\maketitle
\thispagestyle{empty}
\pagestyle{empty}

\begin{abstract}
Motivated by the spread of opinions on different social networks, we study a distributed continuous-time bi-virus model
for a system of groups of individuals.
An in-depth stability
analysis is performed for more general models than have been previously considered, 
for the healthy and epidemic states.
In addition, we investigate
sensitivity properties of some nontrivial equilibria and obtain an impossibility result
for distributed feedback control.
\end{abstract}


\section{Introduction}

The spread of epidemic processes over large populations is an important research area,
and is in fact a widely studied topic in epidemiology \cite{book}. To model such a process,
various epidemic models have been proposed such as susceptible-infected-recovered (SIR), susceptible-exposed-infected-recovered (SEIR), and susceptible-infected-susceptible (SIS) models \cite{sir,seir,yorke}.
Bernoulli developed one of the first known models inspired by the smallpox virus \cite{bernoulli1760essai}.
In this paper we focus on  
continuous time 
SIS models 
\cite{kermack1932contributions,OmicTN09,ahn2013global,pare2015stability}.  Such models consist of a number of agents that are either infected or healthy (susceptible), and cycle back and forth between these two states depending on their current state, connection to infected neighbors, and infection and healing rates. 

The idea of competing SIS virus models has been  investigated in \cite{prakash2012winner,sahneh2014competitive,wei2013competing,santos2015bivirus}. The main motivation for studying these systems has been to understand how competing opinions spread on different social networks. Competing viruses have also been explored for an SIR 
model in \cite{yagan2013conjoining}. In \cite{prakash2012winner}, a homogeneous virus model is studied, that is, one where the infection and healing rates are the same for all agents, with both viruses propagating over the same undirected graph structure. The set of equilibrium points has been determined and sufficient conditions for local stability are given for all equilibria except the coexisting equilibrium. In \cite{sahneh2014competitive}, the equilibria of  a heterogeneous virus model, with undirected graph structures for both viruses assumed, are studied. Co-existence of the epidemic states from both viruses is shown with, however, no stability analysis provided. In \cite{santos2015bivirus}, a sufficient condition for the survival of a single virus is given for a  model, where it is assumed that the virus, homogeneous in the healing rate, propagates over undirected, regular graphs. Note all  previous work has been conducted for undirected graph structures with limited/local stability analysis.

In this paper, we study a distributed continuous-time bi-virus model over {\em directed} graphs.
Consider $n>1$ groups of individuals, labeled $1$ to $n$.
An individual may be infected by individuals in its own as well as nearest-neighbor groups.
Neighbor relationships among the $n$ groups are described
by a directed graph $\bbb{G}$ on $n$ vertices with an arc (or a directed edge) from vertex $j$
to vertex $i$ whenever the individuals in group $i$ can be infected by those in group $j$.
Thus, the neighbor graph $\bbb{G}$ has self-arcs at all $n$ vertices and the directions of arcs
in $\bbb{G}$ represent the directions of the epidemic contagion.
We assume that $\bbb{G}$ is strongly connected.


One contribution of this paper is the analysis
of the equilibria of the bi-virus model over directed graphs (defined in Section \ref{system})
and their stability under appropriate conditions given in Section~\ref{equilibria}.
A second contribution is the development of a sensitivity result for nontrivial equilibria with respect to the infection rate and healing rate, $\delta_i$ and $\beta_{ij}$, in Section~\ref{sense}.
An interesting and surprising impossibility result
for a distributed feedback controller is also provided in Section~\ref{feedback}.
We now begin with some notation and preliminary results. 

\subsection{Notation}

For any positive integer $n$, we use $[n]$ to denote the set $\{1,2,\ldots,n\}$.
We view vectors as column vectors. We use $x'$ to denote the transpose of a vector $x$ and,
similarly, we use $A'$ for the transpose of a matrix $A$.
The $i$th entry of a vector $x$ will be denoted by $x_i$.
The $ij$th entry of a matrix $A$ will be denoted by $a_{ij}$ and, also, by $[A]_{ij}$ when convenient.
We use $\0$ and $\1$ to denote the vectors whose entries all equal $0$ and $1$, respectively,
and $I$ to denote the identity matrix,
while the dimensions of the vectors and matrices are to be understood from the context.
For any vector $x\in\R^n$, we use ${\rm diag}(x)$ to denote the $n\times n$ diagonal matrix
whose $i$th diagonal entry equals $x_i$.
For any two sets $\scr{A}$ and $\scr{B}$,
we use $\scr{A}\setminus \scr{B}$ to denote the set of elements in $\scr{A}$ but not in $\scr{B}$.

For any two real vectors $a,b\in\R^n$, we write $a\geq b$ if
$a_{i}\geq b_{i}$ for all $i\in[n]$,
$a>b$ if $a\geq b$ and $a\neq b$, and $a \gg b$ if $a_{i}> b_{i}$ for all $i\in[n]$.
Similarly, for any two real matrices $A,B\in\R^{m\times n}$, we write $A\geq B$ if
$A_{ij}\geq B_{ij}$ for all $i\in[m]$ and $j\in[n]$,
$A>B$ if $A\geq B$ and $A\neq B$, and $A \gg B$ if $A_{ij}> B_{ij}$ for all $i\in[m]$ and $j\in[n]$.

For a real square matrix $M$, we use $\rho(M)$ to denote its spectral radius 
and $s(M)$ to denote the largest real part among its eigenvalues, i.e.,
\begin{eqnarray*}
\rho(M) &=& \max \left\{|\lambda|\ : \ \lambda\in\sigma(M)\right\}, \\
s(M) &=& \max \left\{{\rm Re}(\lambda)\ : \ \lambda\in\sigma(M)\right\},
\end{eqnarray*}
where $\sigma(M)$ denotes the spectrum of $M$.

The sign function of a real number $x$ is defined as follows:
$$
{\rm sgn}(x) = \left\{ \begin{array}{lll}
  -1 &\mbox{ if \ $x<0$}, \\
  0 &\mbox{ if \ $x=0$}, \\
  1 &\mbox{ if \ $x>0$}.
       \end{array} \right.
$$
Note that for any real number $x\neq 0$,
$$\frac{d|x|}{dx} = {\rm sgn}(x).$$

\subsection{Preliminaries}

For any two nonnegative vectors $a$ and $b$ in $\R^n$,
we say that $a$ and $b$ have the same sign pattern if they have zero entries
and positive entries in the same places, i.e., for all $i\in[n]$,
$a_i=0$ if and only if $b_i=0$, and $a_i>0$ if and only if $b_i>0$.
A square matrix is called {\em irreducible}  if it cannot be permuted to a block upper triangle matrix.

\begin{lemma}
Suppose that $Mx=y$ where $M\in\R^{n\times n}$ is an irreducible nonnegative matrix
and $x,y>\0$ are two vectors in $\R^n$.
If $x$ has at least one zero entry, then $x$ and $y$ cannot have the same sign pattern.
In particular, there exists an index $i\in[n]$ such that $x_i=0$ and $y_i>0$.
\label{pattern}\end{lemma}

{\em Proof:}
Suppose, to the contrary, that for all $i\in[n]$ such that $x_i=0$,  $y_i=0$.
Since $x$ has at least one zero entry and $x>\0$, there exists a proper nonempty subset
$\scr{E}\subset [n]$ such that for any $i\in\scr{E}$, $x_i=0$,
and for any $i\in [n] \setminus \scr{E}$, $x_i>0$. It follows that for any $i\in\scr{E}$, $y_i=0$.
Without loss of generality,
set $\scr{E}=\{1,2,\ldots,m\}$, $1\le m<n$.
Then,
$$Mx = \begin{bmatrix}
A & B \cr
C & D\end{bmatrix}=
\begin{bmatrix}\; \0 \; \cr \bar x\end{bmatrix}
= y = \begin{bmatrix} \; \0 \; \cr \bar y\end{bmatrix},$$
with $A,D\ge0$, $B,C>0$, $\bar x \gg \0$, and $\bar y \ge \0$.
This implies B = 0, which is a contradiction since $M$ is an irreducible matrix.
Therefore, there exists an index $i\in[n]$ such that $x_i=0$ and $y_i>0$,
and thus $x$ and $y$ cannot have the same sign pattern.
\hfill
$\qed$

A real square matrix is called {\em Metzler} if its off-diagonal entries are all nonnegative.
We will use the following important properties of Metzler matrices.

\begin{lemma}
{\rm (Lemma 2.3 in \cite{varga})}
Suppose that $M$ is an irreducible Metzler matrix. Then,
$s(M)$ is a simple eigenvalue of $M$ and there exists a unique (up to scalar multiple) vector $x\gg \0$ such that
$Mx=s(M)x$.
\label{metzler0}\end{lemma}

\begin{lemma}
{\rm (Section 2.1 in \cite{varga})}
Suppose that $M$ is an irreducible Metzler matrix in $\R^{n\times n}$ and $x>\0$ is a vector in $\R^n$.
If $Mx < \lambda x$, then $s(M)<\lambda$.
If $Mx = \lambda x$, then $s(M)=\lambda$. If $Mx > \lambda x$, then $s(M)>\lambda$.
\label{metzler}\end{lemma}

\begin{proposition}
Suppose that $\Lambda$ is a negative diagonal matrix in $\R^{n\times n}$
and $N$ is an irreducible nonnegative matrix in $\R^{n\times n}$.
Let $M=\Lambda+N$.
Then, $s(M)<0$ if and only if $\rho(-\Lambda^{-1}N)<1$,
$s(M)=0$ if and only if $\rho(-\Lambda^{-1}N)=1$,
and $s(M)>0$ if and only if $\rho(-\Lambda^{-1}N)>1$.
\label{ff}
\end{proposition}



{\em Proof:}
Suppose that $\Lambda$ is a negative diagonal matrix in $\R^{n\times n}$
and $N$ is an irreducible nonnegative matrix in $\R^{n\times n}$.
Let $M=\Lambda+N$. By Theorem 3.29 in \cite{varga},
$s(M)<0$ if and only if $\rho(-\Lambda^{-1}N)<1$.
To prove the proposition, it is enough to show that $s(M)=0$ if and only if $\rho(-\Lambda^{-1}N)=1$.

First suppose that $s(M)=0$.
Set $\Lambda_{\varepsilon}=\Lambda - \varepsilon I$ with $\varepsilon>0$.
Let
$M_{\varepsilon} = \Lambda_{\varepsilon} +N = \Lambda - \varepsilon I +N$.
Then, $\lim_{\varepsilon\rightarrow 0^+}\rho(-\Lambda_{\varepsilon}^{-1}N) = \rho(-\Lambda^{-1}N)$.
Since $\varepsilon>0$, $s(M_{\varepsilon})<0$. Then,
$\rho(-\Lambda_{\varepsilon}^{-1}N) < 1$ and, therefore, $\lim_{\varepsilon\rightarrow 0^+}\rho(-\Lambda_{\varepsilon}^{-1}N) \le 1$.
Thus, $\rho(-\Lambda^{-1}N)\le 1$. To prove that $\rho(-\Lambda^{-1}N) = 1$,
suppose that, to the contrary, $\rho(-\Lambda^{-1}N)<1$.
Then, $s(M)<0$, which is a contradiction. Therefore, $\rho(-\Lambda^{-1}N) = 1$.

Now suppose that $\rho(-\Lambda^{-1}N) = 1$. Again set $\Lambda_{\varepsilon}=\Lambda - \varepsilon I$ with $\varepsilon>0$
and $M_{\varepsilon} = \Lambda_{\varepsilon} +N$.
Then, $\lim_{\varepsilon\rightarrow 0^+} s(M_{\varepsilon}) = s(M)$.
Since $\varepsilon>0$, $-\Lambda_{\varepsilon}^{-1}N$ is a nonnegative matrix.
Since $N$ is irreducible and nonnegative, so is $-\Lambda_{\varepsilon}^{-1}N$.
Note that the $i$th diagonal entry of $-\Lambda_{\varepsilon}$ is strictly larger than
the $i$th diagonal entry of $-\Lambda$ since $\varepsilon>0$. Thus,
$-\Lambda^{-1}N > -\Lambda_{\varepsilon}^{-1}N$.
By the Perron-Frobenius Theorem for irreducible nonnegative matrices,
$\rho(-\Lambda_{\varepsilon}^{-1}N) < 1$.
Then,
$s(M_{\varepsilon})<0$ and, thus, $\lim_{\varepsilon\rightarrow 0^+} s(M_{\varepsilon})\le 0$.
Thus, $s(M)\le 0$. To prove that $s(M) = 0$,
suppose that, to the contrary, $s(M)<0$.
Then, $\rho(-\Lambda^{-1}N)<1$, which is a contradiction. Therefore, $s(M)=0$.
\hfill
$\qed$

\begin{lemma}
{\rm (Proposition 2 in \cite{rantzer})}
Suppose that $M$ is an irreducible Metzler matrix such that $s(M)<0$.
Then, there exists a positive diagonal matrix $P$ such that
$M'P+PM$ is negative definite.
\label{less}\end{lemma}

\begin{lemma}
{\rm (Lemma A.1 in \cite{KhanaferAutomatica14})}
Suppose that $M$ is an irreducible Metzler matrix such that $s(M)=0$.
Then, there exists a positive diagonal matrix $P$ such that
$M'P+PM$ is negative semi-definite.
\label{equal}\end{lemma}

\section{The Bi-Virus Model}\label{system}

As noted in the introduction,
we are interested in the following continuous-time distributed model for two competing viruses,
first proposed in a less general form in \cite{prakash2012winner}:
\begin{equation}\label{updates}
\begin{split}
\dot{x}^{1}_i(t) &= - \delta^{1}_i x^{1}_i(t) + (1 - x^{1}_i(t) - x^{2}_i(t))\sum_{j=1}^n \beta^{1}_{ij}  x^{1}_j(t) , \\
\dot{x}^{2}_i(t) &= - \delta^{2}_i x^{2}_i(t) +(1 - x^{2}_i(t) - x^{1}_i(t))\sum_{j=1}^n \beta^{2}_{ij}  x^{2}_j(t) ,    
\end{split}
\end{equation}
where ${x}^{1}_i(t),{x}^{2}_i(t)$ are the probabilities that agent $i$ has virus $1$ or $2$, respectively, each virus has its own non-symmetric infection rates incorporating the nearest-neighbor graph structures  $\beta^{1}_{ij},\beta^{2}_{ij}$, healing rates $\delta^{1}_i,\delta^{2}_i$, and $x^{1}_i(0),x^{2}_i(0),(1-x^{1}_i(0)-x^{2}_i(0))\in[0,1],i\in[n]$. The model can be written in matrix form as
\begin{equation}\label{sys}
\begin{split}
\dot{x}^{1}(t) &= (- D^{1}+B^{1} - X^{1}(t)B^{1} - X^{2}(t)B^{1} )x^{1}(t), \\
\dot{x}^{2}(t) &= (- D^{2}+B^{2} - X^{2}(t)B^{2} - X^{1}(t)B^{2} )x^{2}(t),
\end{split}
\end{equation}
where ${x}^{k}(t)\in [0,1]^n$, $B^{k}$ is the matrix of $\beta^{k}_{ij}$'s, $X^{k}(t)={\rm diag}(x^{k}(t))$, and $D^{k}={\rm diag}(\delta^{k})$, with $k=1,2$ indicating virus $1$ or $2$. Note that if ${x}^{2}(t) = \0$, the above model recovers the single virus model, 
\eq{\dot z_i(t)=-\delta_i z_i(t)+(1-z_i(t))\sum_{j=1}^n \beta_{ij}z_j(t),
\label{update}}
where $z_i(t)$ is the probability that agent $i$ has the single virus,
$\beta_{ij}$ are the infection rates, $\delta_i$ are the healing rates, and $z_i(0)\in[0,1],i\in[n]$,
or in matrix form
\eq{\dot z(t) = \left(-D+B-Z(t)B \right)z(t). \label{single}}
If we further factor the $\beta_{ij}$ into $\beta_i a_{ij}$, where $\beta_i$ is the infection rate of agent $i$ and the $a_{ij}$ defines the connection structure between agents, we recover the standard, single SIS model \cite{OmicTN09}.
We impose the following assumptions on the parameters.

\begin{assumption}
For all $i\in[n]$, we have $\delta^{1}_i,\delta^{2}_i\ge0$. The matrices $B^{1}$ and $B^{2}$ are nonnegative and irreducible.
\label{para}
\end{assumption}
The nonnegativity assumption on the  matrix $B^{k}$ is equivalent to
$\beta^{k}_{ij}\ge 0$ for all $k\in[2]$ and $i,j\in[n]$.
The assumption of an irreducible matrix $B^{k}$ is equivalent to a strongly connected
spreading graph for virus $k$, $k\in[2]$.


\begin{lemma}
Suppose that Assumption \ref{para} holds. Then, $x^{1}_i(t),x^{2}_i(t),x^{1}_i(t)+x^{2}_i(t)\in[0,1]$ for all $i\in[n]$ and $t\ge 0$.
\label{box}
\end{lemma}

{\em Proof:}
Suppose that at some time $\tau$, $x^{1}_i(\tau),x^{2}_i(\tau),x^{1}_i(\tau)+x^{2}_i(\tau)\in[0,1]$ for all $i\in[n]$.
Consider an index $i\in[n]$.
If $x^{1}_i(\tau)=0$, then from \rep{updates} and Assumption \ref{para}, $\dot x^{1}_i(\tau)\ge 0$.
The same holds for $x^{2}_i(\tau)$ and $x^{1}_i(\tau)+x^{2}_i(\tau)$.
If $x^{1}_i(\tau)=1$, then from \rep{updates} and Assumption \ref{para}, $\dot x^{1}_i(\tau)\le 0$.
The same holds for $x^{2}_i(\tau)$ and $x^{1}_i(\tau)+x^{2}_i(\tau)$.
It follows that $x^{1}_i(t),x^{2}_i(t),x^{1}_i(t)+x^{2}_i(t)$ will be in $[0,1]$ for all times $t\ge \tau$.
Since the above arguments hold for all $i\in[n]$,  $x^{1}_i(t),x^{2}_i(t),x^{1}_i(t)+x^{2}_i(t)$ will be in $[0,1]$ for
all $i\in[n]$ and $t\ge \tau$.
Since it is assumed that $x^{1}_i(0),x^{2}_i(0),x^{1}_i(0)+x^{2}_i(0)\in[0,1]$ for all $i\in[n]$,
it follows that $x^{1}_i(t),x^{2}_i(t),x^{1}_i(t)+x^{2}_i(t)\in[0,1]$ for all $i\in[n]$ and $t\ge 0$.
\hfill
$\qed$

Lemma \ref{box} implies that the set
\begin{eqnarray}\label{D}
    \D &=& \{(x^1, x^2) | x_i^1 \ge 0, x_i^2\ge 0, x_i^1+x_i^2 \le 1 \ \forall i\in [n]\} \nonumber\\
    &=&\{(x^1, x^2) \; | \; x^1\ge \0, \; x^2\ge \0, \; x^1+x^2\le \1\}
\end{eqnarray}
is invariant with respect to the system defined by \rep{sys}. Since $x^{1}_i$ and $x^{2}_i$
denote the probability of agent $i$ being infected by virus 1 or 2, respectively,
and $1-x^{1}_i-x^{2}_i$ denotes the probability of agent $i$ being healthy, it is natural to assume that their initial values are in $[0,1]$,
since otherwise the values will lack any physical meaning for the epidemic model considered here.
Therefore, in this paper,
we will focus on the analysis of \rep{sys} only on the domain $\D$, as defined in \rep{D}.

\section{Equilibria and Their Stability} \label{equilibria}

It can be seen that $x^{1}=x^{2}=\0$ is an equilibrium of the system \rep{sys},
which corresponds to the case when  no individual is infected.
We call this trivial equilibrium the {\em healthy state}.
It will be shown that \rep{sys} also admits nonzero equilibria under appropriate assumptions.
We call those nonzero equilibria {\em epidemic states}.
In this section, we study the stability of the healthy state as well as the epidemic states of \rep{sys}.
To state our results, we need the following definition.

\begin{definition}
Consider an autonomous system
\eq{\dot x(t)  = f(x(t)), \label{def}}
where $f: \scr{X}\rightarrow\R^n$ is a locally Lipschitz map from a domain $\scr{X}\subset\R^n$ into $\R^n$.
Let $z$ be an equilibrium of \rep{def} and $\scr{E}\subset\scr{X}$ be a domain containing $z$.
When the equilibrium $z$ is asymptotically stable such that for any $x(0)\in\scr{E}$ we have 
$\lim_{t\rightarrow\infty}x(t) = z$, then $\scr{E}$ is said to be a domain of attraction for $z$.
\end{definition}

\begin{proposition}
Let $z$ be an equilibrium of \rep{def} and $\scr{E}\subset \scr{X}$ be a domain
containing $z$. Let $V:\scr{E}\rightarrow\R$ be a continuously differentiable function
such that $V(z)=0$, $V(x)>0$ in $\scr{E}\setminus \{z\}$, $\dot V(z)=0$,
and $\dot V(x)<0$ in $\scr{E}\setminus \{z\}$. If $\scr{E}$ is an invariant set,
then the equilibrium $z$ is asymptotically stable with domain of attraction $\scr{E}$.
\label{lya}
\end{proposition}
This proposition is a direct consequence of Lyapunov's stability theorem
(see Theorem 4.1 in \cite{khalil}) and the definition of domain of attraction.

\begin{theorem}
Suppose that Assumption \ref{para} holds.
If $s(-D^{1}+B^{1})\leq 0$ and $s(-D^{2}+B^{2})\leq 0$, then
the healthy state is the unique equilibrium of \rep{sys}, which is
asymptotically stable with  domain of attraction $\D$, as defined in \rep{D}.
\label{0global}\end{theorem}
To prove the theorem, we need the following result for the single virus model \rep{single}.

\begin{proposition}
Consider the single-virus model \rep{single}.
Suppose that $\delta_i\ge0$ for all $i\in[n]$, and that the matrix $B$ is nonnegative and irreducible.
If $s(-D+B)\leq 0$, then
$\0$ is asymptotically stable with domain of attraction $[0,1]^n$.
\label{single0}\end{proposition}
This result has been proved in \cite{FallMMNP07,KhanaferAutomatica14} for the case
when $\delta_i > 0$ for all $i\in[n]$. We extend the result by allowing $\delta_i=0$.

{\em Proof:}
We first consider the case when  $s(-D+B)<0$.
Since $(-D+B)$ is an irreducible Metzler matrix, by Lemma \ref{less},
there exists a positive diagonal matrix $P$ such that
$(-D+B)'P+P(-D+B)$ is negative definite.
Consider the Lyapunov function $V(x(t))=x(t)'Px(t)$.
From \rep{single}, when $x(t)\neq \0$,
\begin{eqnarray*}
\dot V(x(t)) &=& 2x(t)'P\left(-D+B-X(t)B\right)x(t) \\
&< & -2 x(t)'PX(t)Bx(t) \\
&\leq & 0.
\end{eqnarray*}
Thus, in this case, $\dot V(x(t))<0$ if $x(t)\neq \0$.
By Lemma \ref{box} and Proposition \ref{lya},
$x=\0$ is asymptotically stable with  domain of attraction $[0,1]^n$.

Next we consider the case when  $s(-D+B)=0$.
Since $(-D+B)$ is an irreducible Metzler matrix, by Lemma \ref{equal},
there exists a positive diagonal matrix $P$ such that
$(-D+B)'P+P(-D+B)$ is negative semi-definite.
Consider the Lyapunov function $V(x(t))=x(t)'Px(t)$.
From \rep{single},
\begin{eqnarray*}
\dot V(x(t)) &=& 2x(t)'P\left(-D+B-X(t)B\right)x(t) \\
&= & x(t)'\left((-D+B)'P+P(-D+B)\right)x(t) - 2x(t)'PX(t)Bx(t) \\
&\leq & 0.
\end{eqnarray*}
We claim that $\dot V(x(t))<0$ if $x(t)\neq \0$.
To establish this claim, we first consider the case when $x(t)\gg \0$.
Since $B$ is nonnegative and irreducible, $Bx(t)\gg \0$.
Since $P$ is a positive diagonal matrix, it follows that $x(t)'PX(t)Bx(t)>0$,
so $\dot V(x(t))<0$.
Next we consider the case when $x(t)>\0$ and $x(t)$ has at least one zero entry.
Since $(-D+B)$ is an irreducible Metzler matrix and $P$ is a positive diagonal matrix,
$(-D+B)'P+P(-D+B)$ is a symmetric irreducible Metzler matrix.
Since $(-D+B)'P+P(-D+B)$ is negative semi-definite,
it follows that $s((-D+B)'P+P(-D+B))=0$. By Lemma \ref{metzler0}, $0$ is a simple eigenvalue of
$(-D+B)'P+P(-D+B)$ and it has a unique (up to scalar multiple) strictly positive eigenvector
corresponding to the eigenvalue $0$. Thus, $x(t)'\left((-D+B)'P+P(-D+B)\right)x(t)<0$ when
$x(t)>\0$ and $x(t)$ has at least one zero entry.
Therefore, $\dot V(x(t))<0$ if $x(t)\neq \0$.
By Lemma \ref{box} and Proposition \ref{lya},
$x=\0$ is asymptotically stable with  domain of attraction $[0,1]^n$.
\hfill
$\qed$

Now we are in a position to prove Theorem \ref{0global}.

{\em Proof of Theorem \ref{0global}:}
To prove the theorem, it is sufficient to show that
both $x^{1}_i(t)$ and $x^{2}_i(t)$ will asymptotically converge to $\0$ as $t\rightarrow\infty$
for any initial condition.

Since $x^{1}_i(t)$ and $x^{2}_i(t)$ are always nonnegative by Lemma \ref{box},
from \rep{updates}, 
\begin{eqnarray*}
\dot{x}^{1}_i(t) &\le& - \delta^{1}_i x^{1}_i(t) + (1 - x^{1}_i(t))\sum_{j=1}^n \beta^{1}_{ij}  x^{1}_j(t) , \\
\dot{x}^{2}_i(t) &\le& - \delta^{2}_i x^{2}_i(t) + (1 - x^{2}_i(t))\sum_{j=1}^n \beta^{2}_{ij}  x^{2}_j(t) ,
\end{eqnarray*}
which imply that each of the trajectories of $x^{1}_i(t)$ and $x^{2}_i(t)$ is bounded above by
a single-virus model. From Assumption \ref{para} and Proposition \ref{single0},
both $x^{1}_i(t)$ and $x^{2}_i(t)$ will asymptotically converge to $\0$ as $t\rightarrow\infty$,
and thus the healthy state is the unique equilibrium of \rep{sys}.
\hfill
$\qed$

For the healthy state, we can show this condition is also necessary.

\begin{theorem}
Suppose that Assumption \ref{para} holds.
Then, the healthy state is the unique equilibrium of \rep{sys} if and only if
$s(-D^{1}+B^{1})\leq 0$ and $s(-D^{2}+B^{2})\leq 0$.
\label{eq0}\end{theorem}
This theorem is a consequence of the following result.

\begin{proposition}
Consider the single-virus model \rep{single}.
Suppose that $\delta_i\ge0$ for all $i\in[n]$, and that the matrix $B$ is nonnegative and irreducible.
If $s(-D+B)> 0$, then \rep{single} has two equilibria, $\0$ and $x^*$ which satisfies $x^*\gg \0$.
\label{nontrivial}
\end{proposition}
This result has been proved in \cite{FallMMNP07} for the case
when $\delta_i > 0$ for all $i\in[n]$. We extend the result by allowing $\delta_i=0$, inspired by the technique used in \cite{FallMMNP07}.
To prove Proposition \ref{nontrivial}, we need the following lemma.

\begin{lemma}
Consider the single-virus model \rep{single}.
Suppose that $\delta_i\ge0$ for all $i\in[n]$, and that the matrix $B$ is nonnegative and irreducible.
If $x^*$ is a nonzero equilibrium of \rep{single}, then $x^*\gg \0$.
\label{0}
\end{lemma}

{\em Proof:}
Suppose that $x^*$ is a nonzero equilibrium of \rep{single}.
By Lemma \ref{box}, it must be true that
$x^*\geq \0$.
To prove $x^*\gg \0$, suppose that, to the contrary, $x^*$ has at least one zero entry.
Without loss of generality, set $x^*_1=0$. Since $x^*$ is an equilibrium of \rep{single}, from \rep{update},
$$-\delta_1x^*_1 + (1-x^*_1)\sum_{j=1}^n \beta_{1j}x^*_j = \sum_{j=1}^n \beta_{1j}x^*_j = 0.$$
It follows that for any $j\in[n]$ such that $\beta_{1j}>0$, $x^*_j=0$.
By repeating this argument, since $B$ is irreducible,
we have $x^*_i=0$ for all $i\in[n]$. But it contradicts the assumption that $x^*$ is nonzero.
Thus, $x^*\gg \0$.
\hfill
$\qed$

{\em Proof of Proposition \ref{nontrivial}:}
It is enough to show that if $s(-D+B)> 0$, there exists a unique strictly positive equilibrium.
We first show that there exists an $x^*\gg \0$ which is an equilibrium of \rep{single}.

Any equilibrium $x^*$ of \rep{single} must satisfy
$$(-D+B)x^*=X^*Bx^*,$$
or equivalently,
$$Bx^*=Dx^*+X^*Bx^*.$$
Let $c>0$ be any positive constant such that
\eq{s(-D+B) -c > 0. \label{c}}
Such a constant $c$ always exists since $s(-D+B)> 0$.
Set
$$\bar D = D+cI.$$
Then,
$$Bx^*=\bar Dx^*+X^*Bx^* - cx^*.$$
From Assumption \ref{para}, $D$ is a nonnegative diagonal matrix.
Thus, $\bar D$ is nonsingular and $\bar D^{-1}$ is also a positive diagonal matrix.
It follows that
\begin{eqnarray*}
\bar D^{-1}Bx^* &=& x^*+\bar D^{-1}X^*Bx^*-c\bar D^{-1}x^* \\
&=& x^*+X^*\bar D^{-1}Bx^*-c\bar D^{-1}x^* \\
&=& x^*+{\rm diag}(\bar D^{-1}Bx^*)x^*-c\bar D^{-1}x^* \\
&=& (I-c\bar D^{-1}+{\rm diag}(\bar D^{-1}Bx^*))x^*.
\end{eqnarray*}
Since $\bar D=D+cI$ and $D$ is nonnegative, it follows that
$I-c\bar D^{-1}$ is a nonnegative diagonal matrix.
In the case when $x^*\gg \0$, since $D^{-1}$ is a positive diagonal matrix and $B$ is an irreducible nonnegative matrix,
it follows that ${\rm diag}(\bar D^{-1}Bx^*)$ is a positive diagonal matrix. Therefore 
 $(I-c\bar D^{-1}+{\rm diag}(\bar D^{-1}Bx^*))$ a positive diagonal matrix and is invertible, which implies
$$x^*=\left(I-c\bar D^{-1}+{\rm diag}(\bar D^{-1}Bx^*)\right)^{-1}\bar D^{-1}Bx^*.$$

Consider the above equation and define a map $f:(0,1]^n\rightarrow[0,1]^n$ given by
$$f(x) = \left(I-c\bar D^{-1}+{\rm diag}(\bar D^{-1}Bx)\right)^{-1}\bar D^{-1}Bx.$$
Note that the $i$th entry of $f(x)$, denoted by $f_i(x)$, is given by
$$f_i(x) = \frac{\left(\bar D^{-1}Bx\right)_i}{1-\frac{c}{c+\delta_i}+\left(\bar D^{-1}Bx\right)_i}.$$
Since $\bar D^{-1}$ and $B$ are both nonnegative,
for any $y \geq z$ in $(0,1]^n$,  $f_i(y)\geq f_i(z)$,
so $f(y)\geq f(z)$.

Since $\bar D^{-1}B$ is an irreducible nonnegative matrix, from the Perron-Frobenius Theorem,
there exists $v\gg \0$ such that
\eq{\bar D^{-1}Bv = rv, \label{r}}
where
$$r=\rho(\bar D^{-1}B).$$
Since $s(-\bar D +B) = s(-D+B)-c $, from \rep{c}, $s(-\bar D +B)>0$.
By Proposition \ref{ff}, it follows that
$r>1$. Then, we can always find an $\varepsilon>0$
such that for each $i\in[n]$,
\eq{\varepsilon v_i \le \frac{r-1}{r}.\label{small}}
From this,
$$1\le \frac{r}{1+\varepsilon rv_i},$$
and thus,
$$\varepsilon v_i \le \frac{\varepsilon rv_i}{1+\varepsilon rv_i}.$$
From \rep{r},
$$\varepsilon v_i \le \frac{\left(\bar D^{-1}B\varepsilon v\right)_i}{1+\left(\bar D^{-1}B\varepsilon v\right)_i}
\le \frac{\left(\bar D^{-1}B\varepsilon v\right)_i}{1-\frac{c}{c+\delta_i}+\left(\bar D^{-1}B\varepsilon v\right)_i},$$
which implies that $\varepsilon v \le f(\varepsilon v)$.
It follows from \rep{small} that $\varepsilon v \ll \1$.
Since we have shown that for any $y \geq z$ in $(0,1]^n$, 
$f(y)\geq f(z)$.
It follows that
$f$ maps the compact set $\scr{C}=\{x\ | \ \varepsilon v \le x \le \mathbf{1} \}$
to itself. By Brouwer's fixed-point theorem, $f$ has a fixed point in $\scr{C}$,
which must be strictly positive.

To prove the proposition, it remains to be shown that the fixed point is unique.
Suppose that $x$ and $y$ are both nonzero equilibria of \rep{single}. From Lemma \ref{0},
it follows that $x,y\gg \0$. Set
$$\kappa = \max_{i\in [n]} \frac{x_i}{y_i}.$$
Then, $x\le \kappa y$, and there exists a $j\in[n]$ for which $x_j = \kappa y_j$.
We claim that $\kappa \le 1$. To establish this claim, suppose that,
to the contrary, $\kappa >1$.
Since $x$ is a fixed point of $f$ and for any $u \geq v$ in $(0,1]^n$, 
$f_j(u)\geq f_j(v)$ for all $j\in[n]$,
it follows that
\begin{eqnarray*}
x_j &=& \frac{\left(\bar D^{-1}Bx\right)_j}{1-\frac{c}{c+\delta_i}+\left(\bar D^{-1}Bx\right)_j} \\
&\le&  \frac{\left(\bar D^{-1}B\kappa y\right)_j}{1-\frac{c}{c+\delta_i}+\left(\bar D^{-1}B\kappa y\right)_j} \\
&=&
\frac{\kappa\left(\bar D^{-1}B y\right)_j}{1-\frac{c}{c+\delta_i}+\kappa\left(\bar D^{-1}B y\right)_j}.
\end{eqnarray*}
By the assumption $\kappa >1$ we have,
$$\frac{\kappa\left(\bar D^{-1}B y\right)_j}{1-\frac{c}{c+\delta_i}+\kappa\left(\bar D^{-1}B y\right)_j}<
\frac{\kappa\left(\bar D^{-1}B y\right)_j}{1-\frac{c}{c+\delta_i}+\left(\bar D^{-1}B y\right)_j}
.$$
Since $y$ is a fixed point of $f$,
$$
\frac{\left(\bar D^{-1}B y\right)_j}{1-\frac{c}{c+\delta_i}+\left(\bar D^{-1}B y\right)_j}
= y_j.$$
Then, it follows that
$$x_j <
\frac{\kappa\left(\bar D^{-1}B y\right)_j}{1-\frac{c}{c+\delta_i}+\left(\bar D^{-1}B y\right)_j} = \kappa y_j =x_j,$$
which is a contradiction. Therefore, $\kappa \le 1$,
which implies that $x\le y$. Using the same arguments,
it also can be shown that $y\le x$. Thus, $x=y$, which establishes
the uniqueness of the positive equilibrium.
\hfill
$\qed$

{\em Proof of Theorem \ref{eq0}:}
It has been shown in Theorem \ref{0global} that if $s(-D^{1}+B^{1})\leq 0$ and $s(-D^{2}+B^{2})\leq 0$,
the healthy state is the unique equilibrium of \rep{sys}.
Thus, to prove the theorem, it is sufficient to show that if either $s(-D^{1}+B^{1})> 0$ or $s(-D^{2}+B^{2})> 0$,
the system \rep{sys} admits an epidemic state.

Without loss of generality, suppose that $s(-D^{1}+B^{1})> 0$. Set $x^{2}=\0$. Then,
the dynamics of $x^{1}$ simplifies to a single-virus system, which admits an epidemic state
by Proposition \ref{nontrivial}. Therefore, in the case when  $s(-D^{1}+B^{1})> 0$,
the system \rep{sys} always admits an equilibrium of the form $(\tilde x^{1}, \0)$ with $\tilde x^{1}\gg \0$.
\hfill
$\qed$


Now we turn to the analysis of epidemic states.

\begin{theorem}
Suppose that Assumption \ref{para} holds.
If $s(-D^{1}+B^{1})> 0$ and $s(-D^{2}+B^{2})\leq 0$, then
\rep{sys} has two equilibria,
the healthy state $(\0,\0)$, which is asymptotically stable with domain of attraction $\{(\0,x^2) | x^2 \in [0,1]^n\}$, and a unique epidemic state of the form $(\tilde x^{1}, \0)$ with $\tilde x^{1}\gg \0$, which is
asymptotically stable with  domain of attraction $\D\setminus\{(\0,x^2) | x^2 \in [0,1]^n\}$, with $\D$ defined in \rep{D}.
\label{eglobal}\end{theorem}
To prove the theorem, we need the following result for the single-virus model \rep{single}.

\begin{proposition}
Consider the single-virus model \rep{single}.
Suppose that $\delta_i\ge0$ for all $i\in[n]$, and that the matrix $B$ is nonnegative and irreducible.
If $s(-D+B)> 0$, then
the epidemic state $x^*$ is asymptotically stable with  domain of attraction $[0,1]^n\setminus\{\0\}$.
\label{2local}\end{proposition}
This result has been proved in \cite{FallMMNP07,KhanaferAutomatica14} for the case
when $\delta_i > 0$ for all $i\in[n]$. We extend the result by allowing $\delta_i=0$.
To prove Proposition \ref{2local}, we need the following lemma.

\begin{lemma}
Consider the single-virus model \rep{single}.
Suppose that $\delta_i\ge0$ for all $i\in[n]$, and that the matrix $B$ is nonnegative and irreducible.
If $x(0)\neq \0$, then there exists a $\tau \ge 0$ such that $x(\tau)\gg \0$.
\label{attract}
\end{lemma}

{\em Proof:}
If $x(0)\gg \0$, then the lemma is true with $\tau=0$.
Suppose that $x(0)> \0$.
Let $\scr{F}(t)$ be the set of all those labels $i\in[n]$ such that
$x_i(t)=0$. Since $x(0)> \0$, the set $\scr{F}(0)$ is nonempty.
In other words, $x_i(t)=0$ for all $i\in\scr{F}(t)$
and $x_i(t)>0$ for all $i\in[n]\setminus\scr{F}(t)$.
Since $x(0)\neq \0$ and the matrix $B$ is irreducible by Assumption \ref{para},
there exists at least one label $j\in\scr{F}(0)$ such that
$j$ has a nonzero neighbor $k$, i.e., $x_j(0)=0$, $x_k(0)>0$, and $\beta_{jk}>0$.
From \rep{update}, it follows that $\dot x_j(0)>0$.
Thus, there must exist a $\tau_1>0$ such that $x_j(\tau_1)>0$
and $x_i(\tau_1)>0$ for all $i\in[n]\setminus\scr{F}(0)$.
This implies that $\scr{F}(\tau_1)$ is a proper subset of $\scr{F}(0)$.
Note that $\scr{F}(0)$ is a finite set. By repeating the above arguments,
there exists a $\tau>0$ such that $\scr{F}(\tau)$ is the empty set,
which implies that $x_i(\tau)>0$ for all $i\in[n]$.
\hfill
$\qed$

{\em Proof of Proposition \ref{2local}:}
Let $y_i(t) = x_i(t)- x^*_i$ for all $i\in[n]$.
Set $y(t)=x(t)-x^*$ and let $Y(t)={\rm diag}(y_i(t))$ and $X^*={\rm diag}(x^*_i)$.
Note that
\eq{ \left(-D+B-X^*B\right)x^* =\0. \label{2balance}}
Then, 
\begin{eqnarray*}
\dot y(t) &=& \left( -D+B-\left(Y(t)+X^*\right)B \right)\left(y(t)+x^*\right)\\
&=& \left( -D+\left(I-X^*\right)B -Y(t)B \right)y(t) - Y(t)Bx^* \\
&=& \left( -D+\left(I-X^*\right)B \right)y(t) - Y(t)Bx(t) \\
&=& \left( -D+\left(I-X^*\right)B - {\rm diag}(Bx(t))\right)y(t).
\end{eqnarray*}
Thus, for all $i\in[n]$,
$$\dot y_i(t) = -\delta_i y_i(t) + (1-x^*_i)\sum_{j=1}^n \beta_{ij}y_j(t) - \left(\sum_{j=1}^n \beta_{ij}x_j(t)\right)y_i(t).$$

By Lemma \ref{nontrivial}, we have $x^*_i >0$ for all $i\in[n]$.
Consider the Lyapunov function
$$V(y(t))=\max_{k\in[n]}\frac{|y_k(t)|}{x^*_k}.$$
Then, $V(y(t))\ge 0$ with equality if and only if $y(t)=\0$ (or equivalently, $x(t)=x^*$).
For any time $t$, without loss of generality, let $m$ be an index in $[n]$ for which
$$\frac{|y_m(t)|}{x^*_m}=V(y(t))=\max_{k\in[n]}\frac{|y_k(t)|}{x^*_k}.$$
Then, for all $i\in[n]$,
$$|y_i(t)| \le V(y(t)) x^*_i.$$
Thus, when $|y_m(t)|>0$ (or equivalently, $V(y(t))>0$), 

\begin{eqnarray} 
\dot V(y(t)) &=& \frac{1}{x^*_m}\cdot \frac{d|y_m(t)|}{dt} \nonumber\\
&= & \frac{1}{x^*_m} {\rm sgn}(y_m(t))\dot y_m(t) \label{kkk}\\
&= & \frac{1}{x^*_m} {\rm sgn}(y_m(t))\Bigg( -\delta_m y_m(t) + (1-x^*_m)\sum_{j=1}^n \beta_{mj}y_j(t) \nonumber  - \Bigg(\sum_{j=1}^n \beta_{mj}x_j(t)\Bigg)y_m(t)\Bigg) \nonumber\\
&=& \frac{1}{x^*_m} \Bigg( -\delta_m |y_m(t)| + (1-x^*_m)\sum_{j=1}^n \beta_{mj}y_j(t){\rm sgn}(y_m(t)) \nonumber  - \Bigg(\sum_{j=1}^n \beta_{mj}x_j(t)\Bigg)|y_m(t)|\Bigg) \nonumber\\
&\le & \frac{1}{x^*_m} \Bigg( -\delta_m |y_m(t)| + (1-x^*_m)\sum_{j=1}^n \beta_{mj}|y_j(t)|\Bigg) \nonumber  - \frac{1}{x^*_m} \Bigg(\sum_{j=1}^n \beta_{mj}x_j(t)\Bigg)|y_m(t)| \nonumber \\
&\le & \frac{V(y(t))}{x^*_m} \Bigg( -\delta_m x^*_m + (1-x^*_m)\sum_{j=1}^n \beta_{mj}x^*_j\Bigg) \nonumber  - \frac{1}{x^*_m} \Bigg(\sum_{j=1}^n \beta_{mj}x_j(t)\Bigg)|y_m(t)| \nonumber \\
&=& - \frac{1}{x^*_m} \Bigg(\sum_{j=1}^n \beta_{mj}x_j(t)\Bigg)|y_m(t)| \label{jjj}\\
&\le & 0.\nonumber
\end{eqnarray}
From \rep{kkk} and the definition of $y_i(t)$, it is straightforward to verify that
$\dot V(y(t)) = 0$ in the case when $x(t)=\0$.
Next we consider the case when $x(t)\gg \0$.
From \rep{jjj}, since the matrix $B$ is irreducible by Assumption \ref{para} and
$x^* \gg \0$ by Lemma \ref{nontrivial}, it can be seen that $\dot V(y(t)) < 0$
if $x(t)$ does not equal $x^*$.

Recall that by Lemma \ref{attract}, as long as $x(0)\neq \0$, there always exists
a finite time $t_0$ at which $x(t_0)\gg \0$.
From Lemma \ref{box} and Proposition \ref{lya}, to prove the theorem,
it remains to be shown that the system \rep{single} is invariant on a compact subset of $(0,1]^n$.
Without loss of generality,
suppose that $x(0)\gg \0$.
Let $\varepsilon$ be a nonnegative real number such that
$$\varepsilon = \max_{k\in[n]}\frac{|x_k(0)-x^*_k|}{x^*_k}.$$
Consider the compact set
$$\scr{B}= \left\{ x \ | \ V(x(t)) \le \varepsilon\right\} \subset \ [0,1]^n.$$
Note that
$$V(x(t))=\max_{k\in[n]}\frac{|x_k(t)-x^*_k|}{x^*_k}.$$
It follows that
for any $z_1\in\scr{B}$ and $z_2\in[0,1]^n\setminus\scr{B}$, 
$V(z_1)< V(z_2)$.
Since we have shown that $\dot V(t) < 0$ if $x(t)\gg \0$ and $x(t)\neq x^*$,
there cannot exist a trajectory from $x(0)$ to any point in $[0,1]^n\setminus\scr{B}$.
Since the above arguments hold for any $x(0)\gg \0$,
such a compact set $\scr{B}$ always exists and
the system \rep{single} is invariant on $\scr{B}\subset (0,1]^n$.
\hfill
$\qed$

We are now in a position to provide a sketch of the proof for Theorem \ref{eglobal}.

{\em Sketched proof of Theorem \ref{eglobal}:}
From the proof of Theorem \ref{0global}, $x^{2}(t)$ will asymptotically converge to $\0$
as $t\rightarrow\infty$ for all initial values $(x^{1}(0),x^{2}(0)) \in \{(\0,x^2) | x^2 \in [0,1]^n\}$.
From \rep{sys},
$$\dot{x}^{1}(t) = (- D^{1}+B^{1} - X^{1}(t)B^{1}) x^{1}(t) - X^{2}(t)B^{1}x^{1}(t).$$
Thus, we can regard the dynamics of $x^{1}(t)$ as an autonomous system
\eq{\dot{x}^{1}(t) = (- D^{1}+B^{1} - X^{1}(t)B^{1}) x^{1}(t),\label{temp}}
with a vanishing perturbation $- X^{2}(t)B^{1}x^{1}(t)$, which converges to $\0$  as $t\rightarrow\infty$.
From Proposition \ref{2local}, the autonomous system \rep{temp} will asymptotically converge to
a unique epidemic state $\tilde x^{1} \gg \0$ for any $x^{1}(0)\in[0,1]^n\setminus\{\0\}$.
It can then be shown that $(x^{1}(t),x^{2}(t))$ will asymptotically converge to
the unique epidemic state $(\tilde x^{1},\0)$ for any $(x^{1}(0),x^{2}(0))\in\D\setminus\{(\0,x^2) | x^2 \in [0,1]^n\}$, with $\D$ defined in \rep{D}.
\hfill
$\qed$

It is clear from the preceding results that as long as one of $s(-D^{k}+B^{k})$, $k\in\{1,2\}$,
is less than or equal to zero, at most one virus will ultimately spread over the network.
A natural question is whether the two viruses can coexist when  $s(-D^{k}+B^{k})$, $k\in\{1,2\}$,
are both larger than zero.
In the following, we will partially answer this question.
We begin with a result regarding non-coexisting equilibria.

Let $(\tilde x^{1},\tilde x^{2})$ be an equilibrium of \rep{sys}.
Here, both $\tilde x^{1}$ and $\tilde x^{2}$ can be $\0$.
Then, 
the Jacobian matrix of the equilibrium, denoted $J(\tilde x^{1},\tilde x^{2})$ is
\begin{equation}\label{jacob}
J(\tilde x^{1},\tilde x^{2}) = 
\begin{bmatrix}
(I - \tilde X^{1} -\tilde X^{2})B^{1} -D^{1} - \tilde B^{1} & -\tilde B^{1} \cr
- \tilde B^{2} & (I - \tilde X^{1} -\tilde X^{2})B^{2} -D^{2} - \tilde B^{2}\end{bmatrix}, 
\end{equation}
where $\tilde B^{i} = {\rm diag}(B^{i}\tilde x^{i})$, $i\in[2]$.

\begin{theorem}
Suppose that Assumption \ref{para} holds.
If $s(-D^{1}+B^{1})> 0$ and $s(-D^{2}+B^{2})>0$, then
\rep{sys} has at least three equilibria,
the healthy state $(\0,\0)$, and two epidemic states of the form $(\tilde x^{1}, \0)$ with $\tilde x^{1}\gg \0$
and $(\0,\tilde x^{2})$ with $\tilde x^{2}\gg \0$.
The healthy state $(\0,\0)$ is  unstable.
\label{noncoexist}\end{theorem}

{\em Proof:}
The existence of the two epidemic states is an immediate consequence of Proposition \ref{nontrivial}.
The healthy state $(\0,\0)$ is always an equilibrium of \eqref{sys}.
Since by \eqref{jacob}  
$$
J(\0,\0) = \begin{bmatrix}
-D^{1}+B^{1} & 0 \cr
0 & -D^{2}+B^{2}\end{bmatrix},
$$
which is unstable as $s(-D^{1}+B^{1})> 0$ and $s(-D^{2}+B^{2})>0$,
the healthy state $(\0,\0)$ is unstable.
\hfill
$\qed$

It turns out that non-coexisting equilibria may not exist even though
$s(-D^{k}+B^{k})$, $k\in\{1,2\}$,
are both larger than zero, as shown in the following special case.

\begin{assumption}
Viruses 1 and 2 spread over the same strongly connected directed graph $\bbb{G}=([n],\scr{E})$, with $\delta^{1}_i=\delta^{1}>0$ and $\delta^{2}_i=\delta^{2}>0$ for all $i\in[n]$,
and $\beta^{1}_{ij}=\beta^{1}>0$ and $\beta^{2}_{ij}=\beta^{2}>0$ for all $(i,j)\in\scr{E}$.
\label{special1}\end{assumption}
Under Assumption \ref{special1}, it should be clear that 
$D^{1}=\delta^{1}I$, $D^{2}=\delta^{2}I$, $B^{1}=\beta^{1}A$, 
and $B^{2}=\beta^{2}A$,
where $A$ is the adjacency matrix of $\bbb{G}$, which is an irreducible Metzler matrix.

\begin{theorem}
Suppose that Assumption \ref{special1} holds.
Then, coexisting equilibria may exist only if
$\frac{\delta^{1}}{\beta^{1}}= \frac{\delta^{2}}{\beta^{2}}$.
\label{homo1}\end{theorem}
This result has been proved in \cite{prakash2012winner} for the case
when $\bbb{G}$ is an undirected graph. We extend the result by allowing $\bbb{G}$ to be directed.
To prove the theorem, we need the following lemma.

\begin{lemma}
Suppose that Assumption \ref{special1} holds.
If $(\tilde x^{1},\tilde x^{2})$ is an equilibrium of \rep{sys}, 
then $\tilde x^{1}+\tilde x^{2}\ll \1$.
\label{nonzero}\end{lemma}

{\em Proof:}
To prove the lemma, suppose that, to the contrary, for some time $\tau \geq 0$, there exists some $i\in[n]$ such that 
$\tilde x^{1}_i(\tau)+\tilde x^{2}_i(\tau) = 1$.
Then, from \rep{updates},  $\dot{\tilde x}^{1}_i(\tau),\dot{\tilde x}^{2}_i(\tau) < 0$. But this contradicts 
the hypothesis that $(\tilde x^{1}_i,\tilde x^{2}_i)$ is an equilibrium. Therefore, $\tilde x^{1}(\tau)+\tilde x^{2}(\tau)\ll \1$.
\hfill
$\qed$

{\em Proof of Theorem \ref{homo1}:}
To prove the theorem, suppose that, to the contrary, there exists an equilibrium 
$(\tilde x^{1},\tilde x^{2})$ such that $\tilde x^{1},\tilde x^{2}>\0$ 
in the case when $\frac{\delta^{1}}{\beta^{1}} \ne \frac{\delta^{2}}{\beta^{2}}$.
From \rep{sys} and Assumption \ref{special1}, 
\begin{eqnarray*}
(I-\tilde X^{1}-\tilde X^{2})A\tilde x^{1} = \frac{\delta^{1}}{\beta^{1}}\tilde x^{1}, \\
(I-\tilde X^{1}-\tilde X^{2})A\tilde x^{2} = \frac{\delta^{2}}{\beta^{2}}\tilde x^{2}.
\end{eqnarray*}
From Lemma \ref{nonzero}, $(I-\tilde X^{1}-\tilde X^{2})$ is a positive diagonal matrix, 
and thus $(I-\tilde X^{1}-\tilde X^{2})A$ is also an irreducible Metzler matrix. 
Since $\tilde x^{1},\tilde x^{2}>\0$, from Proposition \ref{metzler}, 
$s((I-\tilde X^{1}-\tilde X^{2})A)=\frac{\delta^{1}}{\beta^{1}} = \frac{\delta^{2}}{\beta^{2}}$,
which is impossible because of the hypothesis that $\frac{\delta^{1}}{\beta^{1}} \ne \frac{\delta^{2}}{\beta^{2}}$.
Therefore, coexisting equilibria may exist only if
$\frac{\delta^{1}}{\beta^{1}}= \frac{\delta^{2}}{\beta^{2}}$.
\hfill
$\qed$

Without loss of generality we assume $s(A)>\frac{\delta^{1}}{\beta^{1}} > \frac{\delta^{2}}{\beta^{2}}$.

\begin{theorem}
Suppose that Assumption \ref{special1} holds and that
$s(A)>\frac{\delta^{1}}{\beta^{1}} > \frac{\delta^{2}}{\beta^{2}}$.
Then, system \rep{sys} has three equilibria, the healthy state
$(\0,\0)$ which is unstable, $(\tilde x^{1},\0)$ with $\tilde x^{1}\gg \0$ which is unstable,
and $(\0,\tilde x^{2})$ with $\tilde x^{2}\gg \0$ which is locally exponentially stable.
\label{homo2}\end{theorem}
This result has been proved in \cite{prakash2012winner} for the case
when $\bbb{G}$ is an undirected graph. We extend the result by allowing $\bbb{G}$ to be directed with a  proof technique similar to \cite{prakash2012winner}.

{\em Proof:}
From Theorem \ref{homo1}, the system \rep{sys} cannot have any equilibria of the form 
$(\tilde x^{1},\tilde x^{2})$ with $\tilde x^{1},\tilde x^{2}> \0$.
Thus, if $(\tilde x^{1},\tilde x^{2})$ is an equilibrium of \rep{sys}, at least one of 
$\tilde x^{1}$ and $\tilde x^{2}$ equals $\0$. 
It is clear that $(\0,\0)$ is always an equilibrium. 
Suppose that $\tilde x^{1}=\0$ and $\tilde x^{2}> \0$. Then, from Proposition \ref{nontrivial},
it must be true that $\tilde x^{2}\gg \0$ and is unique.
Similarly, in the case when $\tilde x^{1}>\0$ and $\tilde x^{2}= \0$, 
it must be true that $\tilde x^{1}\gg \0$ and is unique.
Thus, the system \rep{sys} has only three equilibria. 

Next we turn to the local stability of the three equilibria. 
Note that from Assumption \ref{special1}, the hypothesis  
$s(A)>\frac{\delta^{1}}{\beta^{1}} > \frac{\delta^{2}}{\beta^{2}}$ implies 
that $s(-D^{1}+B^{1}), s(-D^{2}+B^{2})>0$.
Then, from Theorem \ref{noncoexist}, the healthy state $(\0,\0)$ is  unstable.

From \rep{jacob}, the Jacobian at $(\tilde x^{1},\0)$ equals
$$\begin{bmatrix}
\beta^{1}(I - \tilde X^{1})A -\delta^{1}I - \beta^{1}{\rm diag}(A\tilde x^{1}) & -\beta^{1}{\rm diag}(A\tilde x^{1}) \cr
0 & \beta^{2}(I - \tilde X^{1} )A -\delta^{2}I \end{bmatrix}.
$$
From \rep{sys} and Assumption \ref{special1}, 
$$(I - \tilde X^{1} )A\tilde x^{1} = \frac{\delta^{1}}{\beta^{1}}\tilde x^{1}.$$
It follows from Lemma \ref{nonzero} that $(I - \tilde X^{1} )A$ is an irreducible Metzler matrix.
Then, from Proposition \ref{metzler}, 
$s((I - \tilde X^{1} )A)= \frac{\delta^{1}}{\beta^{1}}$.
Since $s(A)>\frac{\delta^{1}}{\beta^{1}} > \frac{\delta^{2}}{\beta^{2}}$,
it follows that 
\begin{eqnarray*}
s(\beta^{2}(I - \tilde X^{1} )A -\delta^{2}I) &=& \beta^{2}s((I - \tilde X^{1} )A) -\delta^{2} \\
&= & \beta^{2}\left(\frac{\delta^{1}}{\beta^{1}} - \frac{\delta^{2}}{\beta^{2}}\right) \\
&>& 0,
\end{eqnarray*}
which implies that the Jacobian matrix is unstable. Thus, 
the equilibrium $(\tilde x^{1},\0)$ with $\tilde x^{1}\gg \0$ is  unstable.

From \rep{jacob}, the Jacobian at $(\0,\tilde x^{2})$ equals
$$\begin{bmatrix}
\beta^{1}(I - \tilde X^{2} )A -\delta^{1}I  &  0 \cr
-\beta^{2}{\rm diag}(A\tilde x^{2}) & \beta^{2}(I - \tilde X^{2})A -\delta^{2}I - \beta^{2}{\rm diag}(A\tilde x^{2}) \end{bmatrix}.
$$
Using the same arguments as in the previous paragraph, 
 $s(\beta^{1}(I - \tilde X^{2} )A -\delta^{1}I)<0$.
From \rep{sys} and Assumption \ref{special1}, 
$$(I - \tilde X^{2} )A\tilde x^{2} = \frac{\delta^{2}}{\beta^{2}}\tilde x^{2}.$$
Since $\tilde x^{2}\gg \0$ and $A$ is irreducible, 
it must be true that 
$$\left(\beta^{2}(I - \tilde X^{2})A -\delta^{2}I - \beta^{2}{\rm diag}(A\tilde x^{2})\right)\tilde x^{2}<\0.$$
It follows from Lemma \ref{nonzero} that $\beta^{2}(I - \tilde X^{2})A -\delta^{2}I - \beta^{2}{\rm diag}(A\tilde x^{2})$ is an irreducible Metzler matrix.
Then, from Proposition \ref{metzler}, 
$s(\beta^{2}(I - \tilde X^{2})A -\delta^{2}I - \beta^{2}{\rm diag}(A\tilde x^{2}))<0$,
which implies that the Jacobian matrix is stable. Thus,
the equilibrium $(\0,\tilde x^{2})$ with $\tilde x^{2}\gg \0$ is locally exponentially stable.
\hfill
$\qed$

For possible coexisting equilibria, we have the following interesting result. 

\begin{theorem}
Suppose that Assumption \rep{special1} holds and that $s(A)>\frac{\delta^{1}}{\beta^{1}} = \frac{\delta^{2}}{\beta^{2}}$.
If $(\tilde x^{1}, \tilde x^{2})$ with $\tilde x^{1},\tilde x^{2}>\0$ is an equilibrium 
of \rep{sys}, then $\tilde x^{1},\tilde x^{2}\gg\0$ and 
$\tilde x^{1} = \alpha\tilde x^{2}$ for some constant $\alpha>0$.
\label{parallel1}\end{theorem}

{\em Proof:}
From the proof of Theorem \ref{homo1}, 
\begin{align}\label{equil}
\begin{split}
(I-\tilde X^{1}-\tilde X^{2})A\tilde x^{1} = \frac{\delta^{1}}{\beta^{1}}\tilde x^{1}, \\
(I-\tilde X^{1}-\tilde X^{2})A\tilde x^{2} = \frac{\delta^{2}}{\beta^{2}}\tilde x^{2},    
\end{split}
\end{align}
in which $(I-\tilde X^{1}-\tilde X^{2})A$ is an irreducible Metzler matrix. 
From Lemma \ref{metzler0}, it must be true that $\tilde x^{1},\tilde x^{2}\gg\0$ and
$\tilde x^{1} = \alpha\tilde x^{2}$ for some constant $\alpha>0$.
\hfill
$\qed$

\begin{remark}
Note, from \rep{jacob} and \eqref{equil}, it can be verified that 
$$J(\tilde x^{1},\tilde x^{2})
\begin{bmatrix}
\phantom{-}\tilde x^{1}\; \cr
-\tilde x^{1}\; \end{bmatrix}=
0 \times
\begin{bmatrix}
\phantom{-}\tilde x^{1} \;\cr
-\tilde x^{1} \;\end{bmatrix},$$
i.e., the Jacobian matrix has a zero eigenvalue. Therefore nothing can be said about the local stability of the coexisting equilibria. Simulations indicate that, 
depending on the initial condition, the system can arrive at different equilibria of the form $\tilde x^{1} = \alpha\tilde x^{2}$ for different constants $\alpha>0$.
\end{remark}

A similar result can be established for another special case, as specified by the following assumption.

\begin{assumption}
Viruses 1 and 2 spread over the same strongly connected directed graph $\bbb{G}=([n],\scr{E})$, with  $\delta^{1}_i=\delta^{2}_i>0$ for all $i\in[n]$,
and $\beta^{1}_{ij}=\beta^{2}_{ij}$ for all $(i,j)\in\scr{E}$.
\label{special2}\end{assumption}
Under Assumption \ref{special2}, we have
$D^{1}=D^{2}=D$ and $B^{1}=B^{2}=B$,
where $D$ is a positive diagonal matrix and $B$ is an irreducible Metzler matrix.

\begin{theorem}
Suppose that Assumption \rep{special2} holds and that $s(-D+B)>0$.
If $(\tilde x^{1}, \tilde x^{2})$ is an equilibrium
of \rep{sys} with $\tilde x^{1},\tilde x^{2}>\0$, then $\tilde x^{1},\tilde x^{2}\gg\0$, $\tilde x^{1}+\tilde x^{2}$ is unique, and
$\tilde x^{1} = \alpha\tilde x^{2}$ for some constant $\alpha>0$.
\label{parallel2}\end{theorem}

{\em Proof:}
From \rep{sys} and Assumption \ref{special2}, 
\begin{eqnarray*}
\dot x^{1}(t) +\dot x^{2}(t) = \left(-D + B - (X^{1}(t)+X^{2}(t))B \right) (x^{1}(t)+x^{2}(t)).
\end{eqnarray*}
Thus, the dynamics of $x^{1}(t)+x^{2}(t)$ is equivalent to the single-virus model \rep{single}.
From Proposition \ref{nontrivial}, in the case when $s(-D+B)>0$,
$x^{1}(t)+x^{2}(t)$ has a unique nonzero equilibrium in $[0,1]^n$.
Thus, $\tilde x^{1}+\tilde x^{2}$ is unique. 
From \rep{sys},
\begin{eqnarray*}
\dot x^{1}(t) -\dot x^{2}(t) = -D(x^{1}(t)-x^{2}(t)) + (B - (X^{1}(t)+X^{2}(t))B) (x^{1}(t)-x^{2}(t)).
\end{eqnarray*}
Then, 
$$(-D+B - (\tilde X^{1}+\tilde X^{2})B) (\tilde x^{1}-\tilde x^{2}) =\0.$$
Using the same arguments as those in the proof of Lemma \ref{nonzero}, 
it can be shown that $\tilde x^{1}+\tilde x^{2}\ll \1$. Thus,  
$-D+B - (\tilde X^{1}+\tilde X^{2})B$ is an irreducible Metzler matrix. 
By Lemma \ref{metzler}, since 
$$(-D+B - (\tilde X^{1}+\tilde X^{2})B) (\tilde x^{1}+\tilde x^{2}) =\0,$$
and $x^{1}(t)+x^{2}(t)\gg \0$, 
$s(-D+B - (\tilde X^{1}+\tilde X^{2})B)=0$.
From Lemma \ref{metzler0},  either $x^{1}(t)=x^{2}(t)$ or
$x^{1}(t)-x^{2}(t) = \gamma (x^{1}(t)+x^{2}(t))$
for some constant $\gamma>0$. 
In both cases, it must be true that 
$\tilde x^{1} = \alpha\tilde x^{2}$ for some constant $\alpha>0$, 
and thus $\tilde x^{1},\tilde x^{2}\gg\0$.
\hfill
$\qed$

\section{Sensitivity}\label{sense}

We have shown that in the case when $s(- D^{1}+B^{1})>0$ and $s(- D^{2}+B^{2})\le 0$,
the system \rep{sys} has a unique epidemic state of the form $(\tilde x^{1},\0)$ with $\tilde x^{1}\gg\0$, which is stable.
It can be seen that the value of $\tilde x^{1}$ is independent of the matrices $D^{2}$ and $B^{2}$, but
depends on the matrices $D^{1}$ and $B^{1}$,
or equivalently, the parameters $\delta^{1}_i$ and $\beta^{1}_{ij}$.
A natural question is: how does the equilibrium $\tilde x^{1}$ change when the values of
$\delta^{1}_i$ and $\beta^{1}_{ij}$ are perturbed?
The aim of this section is to answer this question.

From the proof of Theorem \ref{eglobal},
the value of $\tilde x^{1}$ equals the unique epidemic state, denoted $x^*$, of the single-virus model \rep{single}
when $s(-D+B)>0$.
Thus, to answer the question just raised, it is equivalent to
study how the equilibrium $x^*$ changes when the values of
$\delta_i$ and $\beta_{ij}$ are perturbed.

For our purposes, we assume in this section that $\delta_i>0$ for all $i\in[n]$.
Then, by Proposition \ref{ff}, $s(-D+B)>0$ if and only if $\rho(D^{-1}B)>1$.

Suppose that $\rho(D^{-1}B)>1$. By Proposition \ref{nontrivial}, the epidemic state $x^*$ is
the unique nonzero equilibrium of \rep{single}, which satisfies the equation
$$(-D+B-X^*B)x^*=\0.$$
Define the mapping $\Phi$ as follows
$$\Phi(x^*,D,B):=(-D+B-X^*B)x^*.$$
Then, the equation $\Phi(x^*,D,B)= \0$ defines an implicit function $g:\R^{n\times n}\times \R^{n\times n}\rightarrow \R^n$
given by
$$x^*=g(D,B).$$

For each pair of matrices $D$ and $B$ for which $\rho(D^{-1}B)>1$,
there must exist a small neighborhood $\scr{B}$ such that for any pair of
matrices $D+\Delta D$ and $B+\Delta B$ in $\scr{B}$,
$$\rho\left((D+\Delta D)^{-1}(B+\Delta B)\right)>1.$$
Here $\Delta D$ is the $n\times n$ diagonal matrix whose $i$th diagonal entry equals $\Delta \delta_i$, which denotes the perturbation
of $\delta_i$, and
$\Delta B$ is the $n\times n$ matrix whose $ij$th entry equals $\Delta \beta_{ij}$, which denotes the perturbation of $\beta_{ij}$.
Let $x^*+\Delta x^*$ denote the new epidemic state resulting from the perturbations. Then,
\begin{eqnarray*}
\left(-D-\Delta D+B+\Delta B - (X^*+\Delta X^*)(B+\Delta B)\right) (x^*+\Delta x^*) =\0,
\end{eqnarray*}
where $\Delta X^*={\rm diag}(\Delta x^*)$.
By ignoring the higher order $\Delta$ terms, it is straightforward to verify that
\begin{eqnarray}
\left(-D+B-X^*B-{\rm diag}(Bx^*)\right)\Delta x^* 
\approx X^*\Delta \delta + (X^*-I)\Delta Bx^*,\label{dd}
\end{eqnarray}
where $\Delta \delta$ is the vector in $\R^n$ whose $i$th entry equals  $\Delta \delta_i$.
First note that
$$\left(-D+B-X^*B-{\rm diag}(Bx^*)\right)x^* = -{\rm diag}(Bx^*)x^*.$$
Since $B$ is an irreducible nonnegative matrix and $x^*\gg \0$,
it follows that ${\rm diag}(Bx^*)$ is a positive diagonal matrix.
Let $c>0$ be any positive constant such that $c$ is strictly smaller than
the minimal diagonal entry of ${\rm diag}(Bx^*)$. Then,
${\rm diag}(Bx^*)> cI$ and thus $-{\rm diag}(Bx^*)x^* < -cx^*$.
Since $(-D+B-X^*B-{\rm diag}(Bx^*))$ is an irreducible Metzler matrix,
by Lemma \ref{metzler0},  $s(-D+B-X^*B-{\rm diag}(Bx^*))<-c<0$,
which implies that $(-D+B-X^*B-{\rm diag}(Bx^*))$ is nonsingular.
Thus, by the Implicit Function Theorem (see, e.g., pages 204-206 in \cite{chiang}),
 the function $x^*=g(D,B)$ is differentiable in the neighborhood $\scr{B}$.
From \rep{dd} we have
\begin{eqnarray*}
\Delta x = \left(-D+B-X^*B-{\rm diag}(Bx^*)\right)^{-1}X^*\Delta \delta +
\left(-D+B-X^*B-{\rm diag}(Bx^*)\right)^{-1}(X^*-I)\Delta Bx^*.
\end{eqnarray*}

To proceed, we need the following lemma.

\begin{lemma}
{\rm (Theorem 2.7 in Chapter 6 of \cite{siam})}
Suppose that $M$ is a nonsingular, irreducible Hurwitz Metzler matrix.
Then, $M^{-1}\ll 0$.
\label{inverse}
\end{lemma}
From this lemma and the preceding discussion, it follows immediately that $(-D+B-X^*B-{\rm diag}(Bx^*))^{-1}$
is a strictly negative matrix.
Since $\0\ll x^* \ll \1$, it follows that
all $x_i^*$ decreases as any $\delta_i$ increases or any $\beta_{ij}$ decreases.
We have proved the following result.

\begin{proposition}
Consider the single-virus model \rep{single}.
Suppose that $\delta_i>0$ for all $i\in[n]$, and that the matrix $B$ is nonnegative and irreducible.
If $s(-D+B)>0$,
then each entry of the epidemic state $x^*$ is a strictly decreasing function of $\delta_i$, $i\in[n]$,
and a strictly increasing function of $\beta_{ij}$, $i,j\in[n]$.
\label{change1}
\end{proposition}

Similarly, we have the following result for the bi-virus system \rep{sys}.

\begin{theorem}
Suppose that $\delta^{1}_i>0$, $\delta^{2}_i\ge 0$ for all $i\in[n]$, and that matrices $B^{1}$ and $B^{2}$ are nonnegative and irreducible.
If $s(-D^{1}+B^{1})>0$ and $s(-D^{2}+B^{2})\le 0$,
then each entry of the epidemic state $\tilde x^{1}$ is a strictly decreasing function of $\delta^{1}_i$, $i\in[n]$,
and a strictly increasing function of $\beta^{1}_{ij}$, $i,j\in[n]$.
\label{change}
\end{theorem}

\section{Distributed Feedback Control}\label{feedback}

In this section, we regard each healing rate as a local control input of each agent $i$.
We begin with the single-virus model \rep{single}.

Suppose that the matrix $B$ is fixed. Let $\delta_i=\sum_{j=1}^n \beta_{ij}$.
Then, the row sums of $D^{-1}B$ all equal $1$.
By the Perron-Frobenius Theorem for irreducible nonnegative matrices (see Theorem 2.7 in \cite{varga}),
$\rho(D^{-1}B)=1$, which is equivalent to $s(-D+B)=0$ because of Proposition \ref{ff}. Thus, by Proposition \ref{single0}, the healthy state $\0$ is asymptotically stable
in this case.
This observation implies that in the case when local control inputs $\delta_i$'s are constant,
there always exist sufficiently large $\delta_i$'s which can stabilize the heathy state.

In the following, we will consider the local control inputs of the form
\eq{\delta_i(t)=k_ix_i(t), \;\;\;\;\; i\in[n],\label{feed}}
where $k_i$ is a feedback gain. Designing the controller as a function of the infection rate $x_i(t)$ is an intuitive approach since if the virus is eradicated, no control should be necessary. In implementation, this could be thought of as a treatment plan for individuals via administration of antidote or alternate treatment techniques.
By~\rep{update}, the system reduces to
\begin{eqnarray*}
\dot x_i(t) = -k_i (x_i(t))^2+(1-x_i(t))\sum_{j=1}^n \beta_{ij}x_j(t), \;\;\;\;\; x_i(0)\in[0,1],\;\;\;\;\; i\in[n].\label{update3}
\end{eqnarray*}
The resulting $n$ state equations can be combined to give
\eq{\dot x(t) = \left(-KX(t)+B-X(t)B \right)x(t), \label{sys2}}
where $K= {\rm diag}([k_1,\dots , k_n])$.
Similar to the original system \rep{single}, we impose the following assumption on the parameters of the new system \rep{sys2}.

\begin{assumption}
For all $i\in[n]$, we have $k_i>0$ and the matrix $B$ is nonnegative and irreducible.
\label{para2}
\end{assumption}


Using the same arguments as in the proof of Lemma \ref{box}, it is straightforward to
verify that the set $[0,1]^n$ is still invariant in the new system \rep{sys2}.
Since both $K$ and $X(t)$ are diagonal matrices, they commute. Then, from \rep{sys2},
\begin{eqnarray*}
\dot x(t) &=& (-X(t)K+B-X(t)B )x(t) \\
&=& (B-X(t)(K+B) )x(t) \\
&=& (-K+(K+B)-X(t)(K+B) )x(t).
\end{eqnarray*}
Thus, the system \rep{sys2} has the same form as the original system \rep{single},
with $D$ and $B$ being replaced by $K$ and $K+B$, respectively.

Note that $K^{-1}(K+B)=I+K^{-1}B$. Since by Assumption \ref{para2}, $K$ is a positive diagonal matrix
and $B$ is an irreducible nonnegative matrix, $K^{-1}$ is a positive diagonal matrix and, thus,
$K^{-1}B$ is an irreducible nonnegative matrix.
By the Perron-Frobenius Theorem for irreducible nonnegative matrices (see Theorem 2.7 in \cite{varga}),
$\rho(K^{-1}B)>0$ and, thus, $\rho(I+K^{-1}B)>1$.
This observation implies, by Proposition \ref{2local},
that the new system \rep{sys2} has a unique nonzero $x^*$ which satisfies
$\0\ll x^* \ll \1$ and is asymptotically stable with  domain of attraction
$[0,1]^n\setminus \{\0\}$.
We are thus led to the following result.

\begin{proposition}
Let Assumption \ref{para2} hold, and let $x(0)> \0$.
Then, for any local control inputs of the form
\rep{feed}, the healthy state $\0$ is not a reachable state of the system \rep{sys2}.
\label{im}
\end{proposition}

Note that instability of the healthy state can also be shown using the Jacobian. However the above shows that the origin is not only an unstable equilibrium but is a repeller, that is, a perturbation in any direction  will drive the state to $x^*\gg \0$.

Now we turn to the bi-virus model \rep{sys}.
We consider the local control inputs of the form
\eq{\delta^{1}_i(t)=k^{1}_ix^{1}_i(t), \;\;\;\;\; \delta^{2}_i(t)=k^{2}_ix^{2}_i(t), \;\;\;\;\; i\in[n],\label{feed2}}
where $k^{1}_i$ and $k^{2}_i$ are feedback gains.
By~\rep{updates}, the system reduces to
\begin{eqnarray*}
\dot{x}^{1}_i(t) &=& - k^{1}_i (x^{1}_i(t))^2 + (1 - x^{1}_i(t) - x^{2}_i(t))\sum_{j=1}^n \beta^{1}_{ij}  x^{1}_j(t) , \\
\dot{x}^{2}_i(t) &=& - k^{2}_i (x^{2}_i(t))^2 +(1 - x^{2}_i(t) - x^{1}_i(t))\sum_{j=1}^n \beta^{2}_{ij}  x^{2}_j(t) ,
\end{eqnarray*}
The above equations can be combined into matrix form:
\begin{equation}\label{sysfb}
\begin{split}
\dot{x}^{1}(t) &= (- K^{1}+(K^{1}+B^{1})-X^{1}(t)(K^{1}+B^{1}))x^{1}(t) - X^{2}(t)B^{1}x^{1}(t), \\
\dot{x}^{2}(t) &= (- K^{2}+(K^{2}+B^{2})-X^{2}(t)(K^{2}+B^{2}))x^{2}(t)  - X^{1}(t)B^{2}x^{2}(t),
\end{split}
\end{equation}
where $K^{1}$ and $K^{2}$ are the $n\times n$ diagonal matrices with the $i$-th diagonal entry equal to $k^{1}_i$
and $k^{2}_i$, respectively.
Similar to the original system \rep{sys}, we impose the following assumption on the parameters of the new system \rep{sysfb}.

\begin{assumption}
For all $i\in[n]$, we have $k^{1}_i,k^{2}_i>0$ and the matrices $B^{1}$ and $B^{2}$ are nonnegative and irreducible.
\label{para3}
\end{assumption}

From the preceding discussion, we have the following.

\begin{theorem}
Let Assumption \ref{para3} hold.
Then, for any local control inputs of the form
\rep{feed2}, the healthy state $(\0,\0)$ is an unstable equilibrium of the system \rep{sys2}.
\label{im2}
\end{theorem}


\section{Conclusion}\label{conclusion}

In this paper we have explored the equilibria of a continuous-time bi-virus model and in so doing, as a by-product we have  improved on the results for the single-virus model. We have provided necessary and sufficient conditions for convergence to the healthy state of \eqref{sys}. We have also provided  results on the epidemic states of \eqref{sys}, including several sufficient conditions for stability and instability, as well as a sensitivity condition. 
We have shown that a  distributed proportional controller of the form $\delta_i(t)=k_ix_i(t)$, can never drive the virus model to the healthy state.  For future work, we will study bi-virus models with time--varying graph structure, similar to \cite{pare2015stability} for the single virus case. We also will analyze the multi-virus case, i.e., more than two competing viruses, for the healthy and epidemic states.

\section{Acknowledgement}

The authors wish to thank Xudong Chen, Daniel Liberzon, and Meiyue Shao (Lawrence Berkeley National
Laboratory) for useful discussions which have contributed to
this work.

\bibliographystyle{unsrt}
\bibliography{bib}

\end{document}